\documentclass[11pt]{amsart}

\usepackage[utf8x]{inputenc} 
\usepackage[a4paper]{geometry} 
\usepackage{enumerate} 
\usepackage{amsmath} 
\usepackage{amssymb}
\usepackage{color,colortbl}
\usepackage[dvipsnames,table]{xcolor}
\usepackage{fancyvrb}   
 
\definecolor{airforceblue}{rgb}{0.36, 0.54, 0.66}
\definecolor{bleudefrance}{rgb}{0.19, 0.55, 0.91}
\definecolor{darkorchid}{rgb}{0.6, 0.2, 0.8}
\definecolor{darkorange}{rgb}{1.0, 0.55, 0.0}
\definecolor{darkspringgreen}{rgb}{0.09, 0.45, 0.27}
\definecolor{commentoutput}{rgb}{0.40, 0.00, 0.0}
\definecolor{output}{rgb}{0.8, 0.0, 0.0}  
\definecolor{circOut}{rgb}{0.4, 1.0, 0.0}

\definecolor{Gray}{gray}{0.9}  
\usepackage[arrow,curve,matrix,tips,frame]{xy}
\usepackage{adjustbox}
\usepackage{multirow}
\usepackage{url}
\usepackage{hyperref}  
                 \hypersetup{ pdfborder={0 0 0}, 
                              colorlinks=true, 
                              linktoc=page,
                              pdfauthor={Giovanni Staglian{\`o}}, 
                              pdftitle={On some families of Gushel-Mukai fourfolds}} 
\usepackage{verbatim}

\theoremstyle{plain} 
\newtheorem{proposition}{Proposition}[section] 
\newtheorem{theorem}[proposition]{Theorem} 
 
\newtheorem{corollary}[proposition]{Corollary} 
 
\theoremstyle{definition}

\newtheorem{example}[proposition]{Example} 
\theoremstyle{remark} 
\newtheorem{remark}[proposition]{Remark} 
                                    
\newcommand{\PP}{{\mathbb{P}}}  

\numberwithin{equation}{section}

\title{On some families of Gushel--Mukai fourfolds}

\address{Dipartimento di Matematica e Informatica, Universit\`a degli Studi di Catania} 
\author[G. Staglian\`o]{Giovanni Staglian\`o}
\email{giovannistagliano@gmail.com}
\subjclass[2010]{14J35, 
                 14J45, 
                 68W30, 
                 14Q10
                } 

\begin{document}

\begin{abstract} 
We give explicit descriptions of some Noether--Lefschetz divisors in the moduli space of Gushel--Mukai fourfolds. As a consequence we obtain that their Kodaira dimension is negative.
\end{abstract}

\maketitle

\section*{Introduction}
In this paper, 
we consider  (complex) 
\emph{ordinary Gushel--Mukai fourfolds}, that is 
quadric hypersurfaces in a transverse intersection 
of the Grassmannian $\mathbb{G}(1,4)\subset\PP^9$ in its Pl\"ucker embedding
with a hyperplane.
The interest in these fourfolds,
especially for rationality questions, is classical (see  e.g.
\cite{Roth1949}).
They also appeared in the  classification 
of 
Fano varieties of coindex $3$ (see \cite{mukai-biregularclassification}).

There are strong similarities between these
 fourfolds 
and cubic fourfolds, \emph{i.e.} smooth cubic hypersurfaces in $\PP^5$.
 For instance, 
as in the case of cubic fourfolds, it is known that 
all Gushel--Mukai fourfolds are unirational. Some rational examples are classical 
and no examples have yet been proved to be irrational.
The similarity with  the case of cubic fourfolds
became more evident thanks to the works \cite{DIM,DK1,DK2,DK3},
where the authors, following  
Hassett's study of cubic fourfolds (see \cite{Hassett,Has00}),
introduced, via Hodge theory and the period map, 
the analogous definition of the \emph{Noether--Lefschetz locus} 
inside the moduli space of Gushel--Mukai fourfolds,
which is a
countably infinite union of irreducible divisors.
They also gave explicit geometric descriptions for the first three components of this locus
(see \cite[Section~7]{DIM} and Example~\ref{eseDIM} below).

The present paper is inspired by the work of Nuer \cite{Nuer},
where  
he gave  explicit descriptions for the first components 
of the Noether--Lefschetz locus in the moduli space of cubic fourfolds,
deducing that they are unirational and hence of negative Kodaira dimension.
Here we do basically the same in the case of Gushel--Mukai fourfolds. 
Our contribution is summarized in Theorem~\ref{mainthm} and Table~\ref{tabella1}, and 
as a consequence we deduce 
that 
 the first ten components 
of the Noether--Lefschetz locus 
in the moduli space of Gushel--Mukai fourfolds 
are uniruled and hence of negative Kodaira dimension (see Corollaries~\ref{cor2} and \ref{corDebarre}).
Notice that in the case of cubic fourfolds, thanks to the work \cite{TaVaA},
there are only a few components of the Noether--Lefschetz locus 
for which we have no information about their Kodaira dimension.

This short paper is organized as follows. In Section~\ref{prelim}
we recall some general facts about Gushel--Mukai fourfolds and del Pezzo fivefolds.
In Section~\ref{sec describ} we present our main results. 
In Section~\ref{sec computations} we introduce a \emph{Macaulay2} package,
which is essential to verify many of our claims.
\subsection*{Acknowledgements}
  The author is grateful to F. Russo and M. Bolognesi for stimulating discussions,
  and to 
  O. Debarre and A. Kuznetsov 
  for pointing out Corollary~\ref{corDebarre}.

\section{Preliminaries}\label{prelim}
In this section, we first recall some general facts about Gushel--Mukai fourfolds 
which were proved in \cite{DIM} (see also \cite{DK1,DK2,DK3}). 
We will also introduce the main 
 tools for our constructions in the next section.
\subsection{Parameter space}
A Gushel--Mukai fourfold $X\subset \PP^8$, GM fourfold for short,
is a smooth Fano fourfold of degree $10$ with $\mathrm{Pic}(X) = \mathbb{Z}\langle \mathcal{O}_X(1)\rangle $ 
and $K_X\in|\mathcal{O}_X(-2)|$.
Equivalently,
$X$ is  a quadratic section of a 
$5$-dimensional linear section 
of the cone in $\PP^{10}$ over the Grassmannian $\mathbb{G}(1,4)\subset \PP^9$ of lines in $\PP^4$ 
(see \cite{mukai-biregularclassification} and also \cite{Gu}).
The group of automorphisms of a GM fourfold is a finite subgroup of $\mathrm{PGL}(9,\mathbb{C})$, 
and
GM fourfolds are parameterized, up to isomorphism, by 
the points of a coarse moduli space $\mathcal M_4^{GM}$ of dimension $24$.
Every $[X]\in \mathcal M_4^{GM}$, outside a locus of codimension $2$, 
can be realized as a quadratic section of a (uniquely determined) del Pezzo fivefold $Y\subset\PP^8$,
that is of a smooth hyperplane section of $\mathbb{G}(1,4)\subset \PP^9$. 
Such fourfolds are called \emph{ordinary} GM fourfolds, 
and they are the only ones in which we are interested in this paper.
Recall also that all smooth hyperplane sections of $\mathbb{G}(1,4)\subset \PP^9$ are projectively equivalent.
 
 \subsection{Noether--Lefschetz locus in \texorpdfstring{$\mathcal M_4^{GM}$}{M4GM}}
 We have a \emph{period map} $\mathfrak{p}:\mathcal M_4^{GM}\to\mathcal{D}$ to a $20$-dimensional quasi-projective variety $\mathcal D$,
 which is dominant with irreducible $4$-dimensional fibers (see \cite[Corollary 6.3]{DK3}). 

For a very general GM fourfold $[X]\in \mathcal M_4^{GM}$, the natural inclusion 
\begin{equation}\label{naturalInclusion}
  H^4(\mathbb G(1,4),\mathbb Z)\cap H^{2,2}(\mathbb G(1,4))\subseteq  H^4(X,\mathbb Z)\cap H^{2,2}(X)
\end{equation}
of middle Hodge groups is an equality.
A GM fourfold $X$ is said to be \emph{special} if the inclusion \eqref{naturalInclusion} is strict. 
This means that the fourfold $X$ contains a surface whose cohomology class ``does not come'' from the Grassmannian $\mathbb{G}(1,4)$.
The set of special GM fourfolds is called the Noether--Lefschetz locus.

Special GM fourfolds 
correspond to a countable union of hypersurfaces $\bigcup_d \mathcal{D}_d\subset \mathcal{D}$, labelled 
by the integers  
$d\geq 10$ with
$d\equiv 0,2$, or $4$ (mod $8$) (see \cite[Lemma~6.1]{DIM}). The hypersurface 
$\mathcal{D}_d\subset\mathcal D$ is 
irreducible if $d\equiv 0$ (mod $4$), and it has two irreducible components $\mathcal D_d'$ and $\mathcal D_d''$
if  $d\equiv 2$ (mod $8$) (see \cite[Corollary~6.3]{DIM});
the same holds true for the hypersurface  $\mathfrak{p}^{-1}(\mathcal{D}_d)\subset\mathcal{M}_4^{GM}$.

\subsection{Discriminants}\label{sec discriminant}
Following \cite[Section 7]{DIM}, suppose that an (ordinary) GM fourfold $X\subset Y=\mathbb{G}(1,4)\cap\PP^8$ 
contains a smooth surface $S$.
Let $a,b$ denote integers such that  $[S] = a\sigma_{3,1}+b\sigma_{2,2}$ 
is the class of $S$ in the Chow ring of $\mathbb{G}(1,4)$.
Then we have that 
$[X]\in\mathfrak{p}^{-1}(\mathcal D_d)$, where $d$ is
the determinant (also called discriminant) 
of the intersection matrix with respect to $(\sigma_{1,1|X}, \sigma_{2|X}-\sigma_{1,1|X}, [S])$,
that 
is 
\begin{equation}\label{discriminant}
d = \det  \begin{pmatrix}
     2&0&b\\
     0&2&a-b\\
     b&a-b&(S)_X^2\end{pmatrix} = 
     4 (S)_X^2-2(b^2+(a-b)^2) .
\end{equation}
The self-intersection of $S$ in $X$ is given by 
\begin{equation}\label{S2X}
(S)_X^2 = 3\,a+4\,b-2\,\deg(S)+4\,g(S)-12\,\chi(\mathcal O_S) + 2\,K_S^2 - 4,
\end{equation}
where $\deg(S)$ and $g(S)$ denote, respectively, the degree and the (sectional) genus of  $S$.

When $d\equiv 2 \ (\mathrm{mod}\ 8)$, we have that 
$[X]\in\mathfrak{p}^{-1}(\mathcal D_d')$ if $a+b$ is even, and $[X]\in\mathfrak{p}^{-1}(\mathcal D_d'')$ 
if $b$ is even.

\subsection{Known explicit descriptions of Noether--Lefschetz divisors}
The first irreducible components of the
Noether–Lefschetz locus in $\mathcal M_{4}^{GM}$
are the following: 
\begin{multline}\label{known}
 \mathfrak{p}^{-1}(\mathcal D_{10}'),\ 
 \mathfrak{p}^{-1}(\mathcal D_{10}''),\ 
 \mathfrak{p}^{-1}(\mathcal D_{12}),\ 
 \mathfrak{p}^{-1}(\mathcal D_{16}),\ 
 \mathfrak{p}^{-1}(\mathcal D_{18}'),\\ 
 \mathfrak{p}^{-1}(\mathcal D_{18}''),\ 
 \mathfrak{p}^{-1}(\mathcal D_{20}),\ 
 \mathfrak{p}^{-1}(\mathcal D_{24}),\ 
 \mathfrak{p}^{-1}(\mathcal D_{26}'),\\ 
 \mathfrak{p}^{-1}(\mathcal D_{26}''),\ 
 \mathfrak{p}^{-1}(\mathcal D_{28}),\ 
 \mathfrak{p}^{-1}(\mathcal D_{32}),\ 
 \mathfrak{p}^{-1}(\mathcal D_{34}'), \ldots
\end{multline}
As far as the author knows, there are 
explicit descriptions only for the first three ones
$\mathfrak{p}^{-1}(\mathcal D_{10}')$ 
$\mathfrak{p}^{-1}(\mathcal D_{10}'')$,
$\mathfrak{p}^{-1}(\mathcal D_{12})$ (see \cite{DIM}),
and recently  for $\mathfrak{p}^{-1}(\mathcal D_{20})$ (see \cite{HoffSta}).
\begin{example}\label{eseDIM} We recall here the explicit descriptions of the loci
$\mathfrak{p}^{-1}(\mathcal D_{10}')$, 
$\mathfrak{p}^{-1}(\mathcal D_{10}'')$,
$\mathfrak{p}^{-1}(\mathcal D_{12})$, which have been given in \cite{DIM}.
\begin{itemize}
\item A two-dimensional linear section of a Schubert variety
$\Sigma_{1,1}\simeq \mathbb{G}(1,3)\subset\mathbb{G}(1,4)$ 
is a quadric surface; its class is $\sigma_{1}^2\cdot \sigma_{1,1} = \sigma_{3,1}+ \sigma_{2,2}$. 
The closure inside $\mathcal M_{4}^{GM}$  of the family 
  of fourfolds containing such a surface coincides with
  $\mathfrak{p}^{-1}(\mathcal{D}_{10}')$. 
\item A two-dimensional linear section of $\mathbb{G}(1,4)$ 
is a quintic del Pezzo surface; its class is $\sigma_{1}^4 = 3\sigma_{3,1} + 2 \sigma_{2,2}$.
The closure inside $\mathcal M_{4}^{GM}$  of the family 
  of fourfolds containing such a surface coincides with
  $\mathfrak{p}^{-1}(\mathcal{D}_{10}'')$.
  \item  A two-dimensional linear section 
  of a Schubert variety
$\Sigma_{2}\subset\mathbb{G}(1,4)$ 
is a cubic scroll; its class is
$\sigma_1^2\cdot \sigma_2 = 2\sigma_{3,1} + \sigma_{2,2}$.
   The closure inside $\mathcal M_{4}^{GM}$  of the family 
  of fourfolds containing such a cubic scroll coincides with
  $\mathfrak{p}^{-1}(\mathcal{D}_{12})$. 
  \end{itemize}
The fourfolds in $\mathfrak{p}^{-1}(\mathcal D_{10}'')$,
as well as some of those in $\mathfrak{p}^{-1}(\mathcal D_{12})$,
have been 
studied in~\cite{Roth1949}.
\end{example}
\begin{example}\label{uniraD20} The description 
of $\mathfrak{p}^{-1}(\mathcal D_{20})$ has been recently achieved in \cite{HoffSta} 
 as the closure 
 inside $\mathcal M_{4}^{GM}$  of the family 
  of fourfolds containing a suitable rational surface  
  of degree $9$ and sectional genus $2$ 
  with class $6\sigma_{3,1} + 3 \sigma_{2,2}$.
  The family of such surfaces in a del Pezzo fivefold 
  is irreducible and unirational of dimension $25$.
  An explicit geometric construction of this family 
  has been provided in \cite{RS3}.
\end{example}
In Section~\ref{sec describ}, we will provide explicit descriptions for all 
the components in the first two rows of \eqref{known}, 
 obtaining some alternative descriptions for those already known.

\subsection{Count of parameters}\label{sec parameter count}
Here we explain how to estimate the dimension of certain families of GM fourfolds 
 just doing some calculations on a specific example. This leads 
 to a computer-implementable algorithm which will be essential for our purposes.
 A similar argument 
 has been given in \cite{Nuer} in the case of families of cubic fourfolds.

Suppose we have a specific example of a smooth surface $S$ 
embedded in a smooth del Pezzo fivefold $Y=\mathbb{G}(1,4)\cap \PP^8\subset\PP^8$, 
and that 
$X\in \mathbb{P}(H^0(\mathcal O_{Y}(2)))$ is a smooth 
quadratic section of $Y$ which contains $S$.
Let $\mathrm{Hilb}^{\chi(\mathcal{O}_S(t))}_{Y}$ be the Hilbert scheme
of subschemes of $Y$ with Hilbert polynomial $\chi(\mathcal{O}_S(t))$, and assume 
that $\mathrm{Hilb}^{\chi(\mathcal{O}_S(t))}_{Y}$ is smooth at $[S]$. 
For instance, this condition is satisfied when we have
\begin{equation}\label{h1}
 h^1(N_{S/Y}) = 0.
\end{equation}
Then there exists a unique irreducible component $\mathcal S$ of 
 $\mathrm{Hilb}^{\chi(\mathcal{O}_S(t))}_{Y}$ 
which contains $[S]$, and the dimension of this component is $h^0(N_{S/Y})$.
Consider the incidence correspondence
\begin{equation*}
\mathcal{Z}={\{([S'], [X'])\;:\; S'\subset X'  \}}\subset \mathcal S\times\mathbb{P}(H^0(\mathcal O_{Y}(2)))
\end{equation*}
and let
\begin{equation}\label{Incidence}
\xymatrix{
 & \mathcal{Z} \ar[dl]_{p_1}\ar[dr]^{p_2}&\\
\mathcal S&&\mathbb{P}(H^0(\mathcal O_{Y}(2)))&}
\end{equation}
be the two 
projections.
We are interested in estimating the codimension of the image 
\begin{equation*}
p_2(\mathcal Z) = { \{ \mbox{quadratic sections of }Y\mbox{ that contain some surface of }\mathcal{S} \} }
\end{equation*}
in the $39$-dimensional projective space 
$\mathbb{P}(H^0(\mathcal O_{Y}(2)))$ of quadratic sections of $Y$, which contains as an open set 
the locus of (ordinary) GM fourfolds contained in $Y$.
By the semicontinuity theorem, the dimension of the fiber $p_1^{-1}([S'])\simeq \PP(H^0(\mathcal I_{S'/Y}(2)))$ of $p_1$ 
at a point $[S']\in \mathcal S$ 
achieves its minimum value on an open set of $\mathcal S$.
Therefore, if $h^0(\mathcal I_{S/Y}(2))$ is minimal,
then we  
deduce that $\mathcal Z$ 
has a unique irreducible component $\mathcal Z^{\circ}$ that dominates $\mathcal S$,
and the dimension of this component is 
\[
 h^0(N_{S/Y}) + h^0(\mathcal I_{S/Y}(2)) - 1.
\]
Note that
$h^0(\mathcal I_{S/Y}(2))$ is automatically minimal if 
\begin{equation}\label{minimality}
 h^1(\mathcal{O}_S(2)) = 0 \mbox{ and } h^0(\mathcal{I}_{S/Y}(2)) = h^0(\mathcal{O}_{Y}(2)) - \chi(\mathcal{O}_S(2)) .
\end{equation}
Now the value of $h^0(N_{S/X})$ is
 the dimension of the tangent space of $\mathrm{Hilb}_X^{\chi(\mathcal{O}_S(t))}$ at
 the special point $[S,X]$. Since $p_2|_{\mathcal Z^{\circ}}^{-1}([X])$ is a closed subscheme of 
 $\mathrm{Hilb}_X^{\chi(\mathcal{O}_S(t))}$, again by 
 semicontinuity, 
we deduce that the dimension of the general fiber of $p_2|_{\mathcal Z^{\circ}}$ is bounded from above 
by $h^0(N_{S/X})$.
It follows that the codimension of $p_2(\mathcal Z^{\circ})$ in $\mathbb{P}(H^0(\mathcal O_{Y}(2)))$ satisfies  
\begin{equation}\label{codim}
 \mathrm{codim}_{\mathbb{P}(H^0(\mathcal O_{Y}(2)))}(p_2(\mathcal Z^{\circ})) \leq 39 - (h^0(N_{S/Y}) + h^0(\mathcal I_{S/Y}(2)) - 1 - h^0(N_{S/X})).
\end{equation}
Notice also that the inequality in \eqref{codim} is actually an equality 
at least in the following case: 
the right-hand side of \eqref{codim} equals~$1$ and the surface
 $S$  satisfies \eqref{discriminant} and \eqref{S2X}
for some value  $d\neq 0$; 
when this happens, 
$p_2(\mathcal{Z}^{\circ})$ 
coincides with one of the (at most two) irreducible components of $\mathfrak{p}^{-1}(\mathcal{D}_d) $,
and
$p_2(\mathcal Z)$ and $p_2(\mathcal Z^{\circ})$ have the same support.

\subsection{Birational representations on \texorpdfstring{$\PP^5$}{P5} of del Pezzo fivefolds}\label{mappeToDP5}
Let $Y= \mathbb{G}(1,4)\cap\PP^8\subset \PP^8$ be a del Pezzo fivefold.
We recall here the two well-known ways 
to obtain a rational parameterization of $Y$.

\subsubsection{} 
Inside $Y$ 
there are two types of planes: a three-dimensional family of 
$\rho$-planes,  having class $\sigma_{2,2}$ in $\mathbb{G}(1,4)$;
and a four-dimensional family of $\sigma$-planes, 
having class $\sigma_{3,1}$ in $\mathbb{G}(1,4)$.
The projection from a $\rho$-plane gives 
a birational map $Y\dashrightarrow\PP^5$, whose inverse 
is defined by the linear system of quadrics through 
a rational normal cubic scroll 
$B=\mathbb{P}(\mathcal{O}_{\PP^1}(2)\oplus \mathcal{O}_{\PP^1}(1))\subset\PP^4\subset\PP^5$.
Conversely, the rational map defined by the quadrics 
through a rational normal cubic scroll $B$ is a birational map
$\psi_{B}:\PP^5\dashrightarrow Y$, whose inverse is the projection from a $\rho$-plane.

The projection of $Y$ from a $\sigma$-plane 
gives 
a dominant map $Y\dashrightarrow Q\subset \PP^5$ onto a smooth quadric hypersurface $Q\subset \PP^5$.

\subsubsection{}
The Hilbert scheme $\mathrm{Hilb}_Y^{2t+1}$ of conics in $Y$ is irreducible  of dimension $10$.
The plane $P$ obtained as the linear span of a general conic $[C]\in \mathrm{Hilb}_Y^{2t+1}$ is not 
contained in $Y$. The projection from $P$ 
gives a birational map $Y\dashrightarrow\PP^5$, 
whose inverse is defined by the linear system of cubics 
through 
a  projection $B$ in $\PP^5$ 
of a rational normal quartic scroll 
$\mathbb{P}(\mathcal{O}_{\PP^1}(2)\oplus \mathcal{O}_{\PP^1}(1)\oplus \mathcal{O}_{\PP^1}(1))\subset\PP^6$.
Conversely, the rational map defined by the cubics 
through a general projection $B$ of a rational normal quartic scroll is a birational map
$\psi_{B}:\PP^5\dashrightarrow Y$, whose inverse is the projection from 
the linear span of a conic in $Y$.

\subsection{A way to construct examples of surfaces in a del Pezzo fivefold}\label{costruire triple}
Let $B\subset\PP^5$ be either a rational normal cubic scroll 
or a rational quartic scroll as considered in Subsection~\ref{mappeToDP5}.
Then we have a birational map $\psi_B:\PP^5\dashrightarrow Y\subset\PP^8$
whose base locus scheme is $B$ and where $Y$ is a del Pezzo fivefold.
A practical way to construct ``\emph{good}'' examples of surfaces in $Y$
consists in taking the image via the map $\psi_B$ 
of some surface $V\subset\PP^5$ cutting $B$ along some curve $C$.

\begin{example}
Let us consider, for instance, the case in which
we want to construct a triple $(B,V,C)$ such that $B\cap V=C$ and where
$B\subset\PP^4\subset\PP^5$ is a rational normal cubic scroll,
$V\subset\PP^5$ is a quintic del Pezzo surface,
and $C\subset\PP^3\subset\PP^5$ a twisted cubic curve.
Then a possible strategy is the following: 
\begin{itemize}
 \item take a twisted cubic curve $C\subset V$ (for example,
 since $V$ can be realized as the image of $\PP^2$ via the linear system of cubic curves with $4$ simple base points,
 we get $C\subset V$ by taking the image of a conic passing through $3$ of the base points);
 \item take a twisted cubic curve $C'\subset B$ (for example, by 
 intersecting $B$ with a general hyperplane);
 \item find a general automorphism $\eta:\PP^5\stackrel{\simeq}{\longrightarrow}\PP^5$ 
        such that $\eta(C')=C$ (for example, 
        by taking standard parametrizations $\nu':\PP^1\to C'\subset\PP^3$ and $\nu:\PP^1\to C\subset\PP^3$,
        one sees that the composition $\nu\circ\nu'^{-1}:C'\to C$ is defined by linear forms 
        and can be extended to an automorphism of $\PP^3\subset\PP^5$);
\item replace $B$ by $\eta(B)$ and we are done.        
\end{itemize}
Now one can verify that the restriction of $\psi_B$ 
induces an isomorphism between $V$ and 
a smooth surface $S=\overline{\psi_B(V)}\subset Y\subset\PP^8$ 
cut out by $19$ quadric hypersurfaces in $\PP^8$
 (see also Remark~\ref{esempi non generali} and the $17$\emph{th} cases in Tables~\ref{tabella1} and \ref{tabella2}).
\end{example}

\begin{example}
 As another example, 
 we construct a triple $(B,V,C)$ such that $B\cap V=C$ and where 
 $B\subset\PP^4\subset\PP^5$ is again a rational normal 
 cubic scroll, while
 $V\subset\PP^5$ is a Veronese surface,
and $C\subset\PP^4\subset\PP^5$ is a rational normal quartic curve.
 Our strategy is the following: 
\begin{itemize}
 \item take a rational normal quartic curve $C\subset V$ (for example, by 
 intersecting $V$ with a general hyperplane);
 \item take a rational normal quartic curve $C'\subset B$ (for example, 
 since $B$ can be realized as the image of $\PP^2$ 
 via the linear system of conics with one base point,
 we get $C'\subset B$ by taking the image of a general conic);
 \item find a general automorphism $\eta:\PP^5\stackrel{\simeq}{\longrightarrow}\PP^5$ 
        such that $\eta(C')=C$ (for example, as before, we can take  
        standard parametrizations $\nu':\PP^1\to C'\subset\PP^4$ and $\nu:\PP^1\to C\subset\PP^4$,
        and then extend $\nu\circ\nu'^{-1}:C'\to C$ to 
        an automorphism of $\PP^4\subset\PP^5$);
\item replace $B$ by $\eta(B)$ and we are done.        
\end{itemize}
One verifies that the restriction of $\psi_B$ 
induces an isomorphism between $V$ and 
another Veronese surface $S=\overline{\psi_B(V)}\subset Y\subset\PP^8$ 
 (see also the $20$\emph{th} cases in Tables~\ref{tabella1} and \ref{tabella2}).
\end{example}

\section{New descriptions of families of GM fourfolds}\label{sec describ}
To construct  surfaces $S$ embedded in a del Pezzo fivefold $Y\subset\PP^8$
we proceed as described in Subsection~\ref{costruire triple}, that is,
we take  the images of  suitable surfaces $V\subset \PP^5$ intersecting the base locus $B$ 
of the map along suitable curves~$C$ via one of the two birational maps $\psi_{B}:\PP^5\dashrightarrow Y\subset\PP^8$ 
described in Subsection~\ref{mappeToDP5}.

More specifically, if we take the triple 
$(B,V,C)$ as in one of the $21$ cases listed in Table~\ref{tabella2}, then the image $S = \overline{\psi_{B}(V)}\subset\PP^8$ 
of the surface $V$ via the map $\psi_{B}:\PP^5\dashrightarrow Y\subset\PP^8$ 
is a smooth irreducible surface cut out scheme-theoretically by hypersurfaces of degree at most two,
contained in smooth quadratic sections of $Y$,
satisfying \eqref{h1} and \eqref{minimality},
and having all the invariants as reported in Table~\ref{tabella1}.

Everything can be verified from the explicit equations of the reducible scheme $B\cup V\subset\PP^5$
using softwares as \emph{Macaulay2} \cite{macaulay2} and \emph{Singular} \cite{singular}.
We give some more details in Section~\ref{sec computations},
by introducing a computational package \cite{SpecialFanoFourfoldsSource} 
that can also produce random reducible schemes $B\cup V$  of the indicated type.
As our results heavily depend on the existence of examples, we also provide
selected examples, one for each case,
on which we have verified all our claims.

Taking into account these examples and 
the facts mentioned in Section~\ref{prelim}, 
we can now formulate the main result of this paper.
\begin{theorem}\label{mainthm}
 Table~\ref{tabella1} contains explicit descriptions 
 for the first $9$ components 
 \begin{equation}\label{comp}
 \mathfrak{p}^{-1}(
 \mathcal D_{10}'
 \cup\mathcal D_{10}'' 
 \cup\mathcal D_{12} 
 \cup \mathcal D_{16}
 \cup\mathcal D_{18}' 
 \cup\mathcal D_{18}'' 
 \cup\mathcal D_{20} 
 \cup\mathcal D_{24} 
 \cup\mathcal D_{26}')
\end{equation}
 of the Noether--Lefschetz locus
 in the moduli space of GM fourfolds.
\end{theorem}
We put in evidence the following immediate consequences.
\begin{corollary}\label{cor1}
 Each irreducible component of 
 \eqref{comp}
can be described as the closure 
  inside $\mathcal M_{4}^{GM}$  of the family 
  of fourfolds containing 
 some  surface 
 of a family of \emph{smooth rational} surfaces in a del Pezzo fivefold. 
\end{corollary}
Since the general member $[S]$  of each of our families of surfaces 
in a del Pezzo fivefold $Y$ satisfies $h^0(\mathcal I_{S / Y}(2)) > 1$,
we deduce that
the general fiber of the left projection 
in \eqref{Incidence} is a positive-dimensional projective space.
In particular, we have the following:
\begin{corollary}\label{cor2}
 Each irreducible component of 
 \eqref{comp}
 is uniruled, and hence of negative Kodaira dimension. 
 (The  1st, 2nd, 3rd, and 7th
 are unirational, see Ex.~\ref{eseDIM} and~\ref{uniraD20}.)
\end{corollary}
\begin{remark}\label{remDebarre}
O. Debarre and A. Kuznetsov pointed out to the author that {the divisors 
$\mathfrak{p}^{-1}(\mathcal D_{d}')$ and $\mathfrak{p}^{-1}(\mathcal D_{d}'')$, 
with $d\equiv 2$ (mod $8$), are birationally isomorphic.} 
This is because the moduli space $\mathcal M_4^{GM}$ 
is a fiber space over the moduli space $LG^0$ 
of Lagrangians $A\subset \bigwedge^3 V_6$ 
with no decomposable vectors, with fiber over a point $[A]$ 
a dense open subset of the quotient 
of the EPW sextic $Y_A^\perp\subset \mathbb{P}(V_6^\vee)$ by its automorphism group 
(\cite[Corollary 6.3]{DK3} and \cite[Proposition 6.9]{DK1}). 
One checks, using the many available results on automorphisms of double EPW sextics, 
that at a general point $[A]$ of any divisor $D\subset LG^0$ 
this automorphism group is trivial. Consider the duality involution 
$\sigma\colon A\mapsto A^\bot$ on $LG^0$. 
Then, the EPW sextic hypersurfaces $Y_A$ and $Y_A^\perp$ are projectively dual 
(hence birationally isomorphic) by \cite[Proposition 6.3]{DK1}, 
hence the divisors $\mathfrak{p}^{-1}(D)$ and $\mathfrak{p}^{-1}(\sigma(D))$ of  $\mathcal M_4^{GM}$ 
are birationally isomorphic. This implies what we need since $\mathcal D_{d}''=\sigma(\mathcal D_{d}')$.
\end{remark}
From the previous remark and Corollary~\ref{cor2}, we immediately deduce the following:
\begin{corollary}\label{corDebarre}
 The divisor $\mathfrak{p}^{-1}(\mathcal D_{26}'')$ 
 is uniruled, and hence of negative Kodaira dimension. 
\end{corollary}

\begin{remark}\label{esempi non generali}
 All our descriptions of loci in $\mathcal M_{4}^{GM}$,
 not only those of codimension one,
are explicit in the sense that each of them is given as 
the closure of the family of fourfolds $X\subset Y$ 
containing some smooth  surface $S\subset Y$ 
which varies on a specific irreducible component $\mathcal S$ of the Hilbert scheme of subschemes of 
a del Pezzo fivefold $Y$.
The existence of the family $\mathcal S$ 
is deduced from the existence of a single example $[S]\in \mathcal S$,
which can be constructed as indicated in Table~\ref{tabella2}. 
Nevertheless, by counting parameters, we see that our constructions 
yield only special examples that  cannot describe
a dense Zariski-open set of $\mathcal S$.
This can be deduced from the last column of Table~\ref{tabella2}.
 For instance, in the $17$\emph{th} case, 
 the family $\mathcal S$ of surfaces $S\subset Y$ of degree 
 $9$ and sectional genus $2$ has dimension $25$. But 
  the family of subschemes of $\PP^5$ which are unions of a del Pezzo surface $V$
and a cubic scroll $B$ intersecting along a twisted cubic curve $C$ has dimension $44$.
So, by sending the surfaces $V$ via the maps $\PP^5\dashrightarrow Y$ defined by the quadrics through $B$,
we get a family of surfaces $S\subset Y$ of 
dimension just $44 - \dim(\mathrm{Aut}(\PP^5)) + \dim(\mathrm{Aut}(Y)) = 44 - 35 + 15 = 24$.
\end{remark}

\begin{remark}
 We are aware of several other examples besides those in Table~\ref{tabella1},
 which can be constructed using the same methods. However, 
 we have found no
 descriptions for the whole divisor 
  $\mathfrak{p}^{-1}(\mathcal D_{26}'')$ and $\mathfrak{p}^{-1}(\mathcal D_d)$ with $d>26$.
  See \cite{rationalSta} for some  sub-loci of $\mathfrak{p}^{-1}(\mathcal D_{26}'')$.
\end{remark}

\subsection{Alternative constructions of surfaces in del Pezzo fivefolds}
Some of our examples of smooth surfaces $S\subset Y=\mathbb{G}(1,4)\cap\PP^8\subset\mathbb{G}(1,4)$ 
  admit an easy interpretation in terms of families of lines.
  We provide two of these examples below.
\begin{example}
Let $X\subset\PP^4$ be a smooth cubic hypersurface. It is classically well-known that $X$ is unirational. Indeed,
let $L\subset X$ be a general line, and consider 
the family $\mathcal{F}_L\subset\mathbb{G}(1,4)$ of lines in $\PP^4$ which are tangent to $X$ at some point of $L$.
Then we get a $2:1$ dominant rational map $\mathcal{F}_L\dashrightarrow X$ 
by sending a line of $\mathcal{F}_L$ to its residual intersection point with $X$.
It follows that $X$ is unirational from the fact that $\mathcal{F}_L$ 
is rational. Actually, $\mathcal{F}_L$ is a cone of vertex $[L]$ 
over a smooth rational normal quartic scroll.
By taking smooth hyperplane sections of $\mathcal{F}_L$, 
we obtain  surfaces in a del  Pezzo fivefold as in case~2 of Table~\ref{tabella1}.
\end{example}
\begin{example}
As shown in \cite{ottavianiScrolls},
 there are only four types of smooth
 threefolds in $\PP^5$ that are scrolls over a surface.
 One of them is a scroll of degree $9$ over a K3 surface of genus $8$.
 Taking a general projection  of such a scroll,
 we get a hypersurface of degree $9$ in $\PP^4$.
 The lines contained in this hypersurface 
 are parameterized by the points of a smooth minimal K3 surface 
 of genus $8$ in $\mathbb{G}(1,4)\cap\PP^8\subset\mathbb{G}(1,4)$,
 which gives an example as in case~3 of Table~\ref{tabella1} (see also \cite{ProcRS}).
\end{example}

\subsection{An application to cubic fourfolds} 
Recall that 
in the $20$-dimensional moduli space $\mathcal C$ of cubic fourfolds,
the Noether--Lefschetz locus consists of a countable infinite union of 
irreducible divisors $\mathcal C_d$, where 
the integer $d$, called the discriminant of the cubic fourfold,
runs over all integers $d\geq 8$ with $d\equiv 0$ or $2$ (mod $6$); 
see \cite{Hassett,Has00,Levico} for details.

If a cubic fourfold $[X]\in \mathcal C$ contains 
an irreducible surface $S$
of degree $\deg(S)$ and sectional genus $g(S)$,
which has smooth normalization and only a finite number $\delta$ of nodes as singularities,
then it is well-known that $[X]\in\mathcal C_d$, where
\begin{equation}\label{formulaDiscr1}
d=3\, (S)_X^2-\deg(S)^2 ,
\end{equation}
and $(S)_X$, the self-intersection of $S$ on $X$,
is given by
\begin{equation}\label{formulaDiscr2}
(S)_X^2 
= 
3\,\deg(S) + 6\,g(S) - 12\,\chi(\mathcal O_S) + 2\,K_S^2 + 2\,\delta - 6 .
\end{equation}

Explicit geometric descriptions of the divisors $\mathcal C_d$
are known only up to $d = 48$: 
the cases $d\leq 44$ with $d\neq 42$ are provided in \cite{Nuer};
the case $d=42$ in \cite{Lai,FarkasVerraC42,RS3};
and the case $d=48$ in \cite{RS3}.
Moreover, as far as the author knows, there are 
no explicit examples of cubic fourfolds known to have discriminant $d > 48$.

We point out that the methods previously used can yield such examples.
We just construct \emph{good surfaces} in a del Pezzo fivefold $Y\subset\PP^8$,
and then we project them onto $\PP^5$ from a general plane $P\subset Y$ of one of the two types of planes in $Y$.
This idea has already been applied in \cite{RS3} 
to get an explicit description of $\mathcal C_{42}$ and deduce 
 the rationality of its general member.
Below is an example among many others that can be constructed in a similar way.
\begin{example}[Some cubic fourfolds of discriminant $60$]
By taking the image of a Veronese surface $V\subset\PP^5$ 
intersecting a cubic scroll $B$ 
along a conic $C$ via the birational map 
$\psi_{B}:\PP^5\dashrightarrow Y\subset\PP^8$ defined by the quadrics through $B$,
we obtain a surface $S\subset Y\subset\PP^8$ 
which is a nodal projection 
of the $3$-uple Veronese embedding  $\nu_3(\PP^2)\subset\PP^9$.
The projection of $S$ from a general $\rho$-plane 
yields a surface in $\PP^5$ of degree $9$ and sectional genus $1$,
cut out by $8$ cubics, having $7$ nodes as the only singularities,
and 
$\nu_3(\PP^2)$ as its normalization.
A general cubic fourfold through such a surface is smooth, and formulas 
\eqref{formulaDiscr1} and \eqref{formulaDiscr2}
tell us that its discriminant equals $60$.
A count of parameters shows that 
the family of these  fourfolds gives
a locus of codimension $3$ in $\mathcal C_{60}$.
\end{example}

\section{Computations: the \emph{Macaulay2} package \emph{SpecialFanoFourfolds}}\label{sec computations}
In this section, we briefly illustrate 
some features of the \emph{Macaulay2} software package 
\href{https://faculty.math.illinois.edu/Macaulay2/doc/Macaulay2/share/doc/Macaulay2/SpecialFanoFourfolds/html/index.html}{\texttt{SpecialFanoFourfolds}} \cite{SpecialFanoFourfoldsSource}, 
with which the data in Table~\ref{tabella1} can be verified.
We refer to its documentation for further details.

This package provides two classes of objects, named \href{https://faculty.math.illinois.edu/Macaulay2/doc/Macaulay2/share/doc/Macaulay2/SpecialFanoFourfolds/html/___Special__Cubic__Fourfold.html}{\texttt{SpecialCubicFourfold}}
and \href{https://faculty.math.illinois.edu/Macaulay2/doc/Macaulay2/share/doc/Macaulay2/SpecialFanoFourfolds/html/___Special__Gushel__Mukai__Fourfold.html}{\texttt{SpecialGushelMukaiFourfold}},
which represent pairs $(S,X)$, with $S$  a surface contained in  a fourfold $X$;
in the former case $X$ is a cubic fourfold, and in the latter case it is a GM fourfold.
An object can be created by passing 
the equations of $S$ and $X$ (or just of $S$ if we want a random $X$ containing $S$).
We now mention some of the functions available, which all take as input just such an object.
\begin{itemize}
\item The function \href{https://faculty.math.illinois.edu/Macaulay2/doc/Macaulay2/share/doc/Macaulay2/SpecialFanoFourfolds/html/_discriminant_lp__Special__Gushel__Mukai__Fourfold_rp.html}{\texttt{discriminant}}, as the name suggests, calculates  
  the discriminant of the fourfold. It applies the 
  formulas \eqref{discriminant}, \eqref{S2X}, \eqref{formulaDiscr1}, \eqref{formulaDiscr2}.
  The 
  Euler-Poincar\'{e} characteristic of the surface, and hence the value of $K_S^2$,
  is calculated using tools from the packages \texttt{Cremona} \cite{packageCremona} 
and \texttt{CharacteristicClasses} \cite{Jost2015}. In the case of GM fourfolds,
the class $a \sigma_{3,1}+b\sigma_{2,2}$ of the surface 
is also obtained automatically through the calculation of an embedding 
into the Grassmannian $\mathbb{G}(1,4)\subset\PP^9$,
 and then applying some tools from the package \texttt{Resultants} \cite{packageResultants}.
\item The function \href{https://faculty.math.illinois.edu/Macaulay2/doc/Macaulay2/share/doc/Macaulay2/SpecialFanoFourfolds/html/_parameter__Count.html}{\texttt{parameterCount}}
automates the count of parameters explained in Subsection~\ref{sec parameter count},
checking in particular 
\eqref{h1} and \eqref{minimality}. 
It gives the needed information on the  
cohomology of the normal bundles. 
We refer to \cite{GraysonBook,Stillman2006ComputingWS} for details on computations with sheaves and sheaf cohomology.
\item  The function \href{https://faculty.math.illinois.edu/Macaulay2/doc/Macaulay2/share/doc/Macaulay2/SpecialFanoFourfolds/html/_detect__Congruence_lp__Special__Gushel__Mukai__Fourfold_cm__Z__Z_rp.html}{\texttt{detectCongruence}} may find  eventual ``\emph{congruences}''
of the surface, which were introduced in \cite{RS1} (see also \cite{RS3}).
In the case
of a special GM fourfold $(S,X)$ contained in a del Pezzo fivefold $Y$, 
the function
counts for each $e\geq1$ the number 
of curves of degree $e$ passing through a general point of $Y$ and that are contained in $Y$ 
and  $(2e-1)$-secants to the surface $S\subset Y$.
If this number is $1$ for some $e$, then 
a \emph{function} is returned, which takes as input  a point $p\in Y$
and returns  the unique $(2e-1)$-secant curve $C\subset Y$ of degree $e$ passing through~$p$.
\item 
Furthermore, after loading the package,
the tools provided by the package \href{https://faculty.math.illinois.edu/Macaulay2/doc/Macaulay2/share/doc/Macaulay2/MultiprojectiveVarieties/html/index.html}{\emph{MultiprojectiveVarieties}} \cite{2021package} are also available. 
In particular, we have two functions that help to construct triples $(B,V,C)$
as the ones considered in Table~\ref{tabella2} (see also Subsection~\ref{costruire triple}):
\begin{itemize}
    \item the function \href{https://faculty.math.illinois.edu/Macaulay2/doc/Macaulay2/share/doc/Macaulay2/MultiprojectiveVarieties/html/_parametrize_lp__Multiprojective__Variety_rp.html}{\texttt{parametrize}} tries 
    to get a rational parametrization of a projective variety (in particular, it works with rational normal curves);
    \item the function \href{https://faculty.math.illinois.edu/Macaulay2/doc/Macaulay2/share/doc/Macaulay2/MultiprojectiveVarieties/html/___Embedded__Projective__Variety_sp_eq_eq_eq_gt_sp__Embedded__Projective__Variety.html}{\texttt{===>}},
    as an application of the function \texttt{parametrize}, tries to find a linear isomorphism 
    between two projective varieties (in particular, it works with rational normal curves of the same degree).
    This function is also used internally (see the function \href{https://faculty.math.illinois.edu/Macaulay2/doc/Macaulay2/share/doc/Macaulay2/SpecialFanoFourfolds/html/_to__Grass.html}{\texttt{toGrass}}) to calculate the embedding of a del Pezzo fivefold into $\mathbb{G}(1,4)\subset\PP^9$.
\end{itemize}
\end{itemize}
We now show how to apply these functions in a specific example. 
For the convenience of the user,
the package also includes a function called \href{https://faculty.math.illinois.edu/Macaulay2/doc/Macaulay2/share/doc/Macaulay2/SpecialFanoFourfolds/html/___G__Mtables.html}{\texttt{GMtables}},
which takes as input an integer $i$ between $1$ and $21$ 
together 
with a coefficient ring $K$, 
and it returns a random triple $(B,V,C)$ defined over $\PP^5_{K}$ 
corresponding to the $i$-\emph{th} case of Table~\ref{tabella2}. 
In the following code, we choose  $i=18$ and $K=\mathbb{Z}/65521$.
(See also \cite{SchreyerAided,Nuer} for the philosophy behind these calculations over finite fields,
although we work here over a finite field only for speed reasons.)
{\footnotesize
\begin{Verbatim}[commandchars=&\[\]]
&colore[darkorange][$ M2 --no-preload]
&colore[output][Macaulay2, version 1.18]
&colore[darkorange][i1 :] &colore[airforceblue][needsPackage] "&colore[bleudefrance][SpecialFanoFourfolds]";
&colore[darkorange][i2 :] (B,V,C) = &colore[bleudefrance][GMtables](18,&colore[darkspringgreen][ZZ]/65521);
\end{Verbatim}
} \noindent 
Now we take the birational map 
$\psi:\PP^5\dashrightarrow Y\subset\PP^8$ defined by the cubics through $B$,
and then we construct 
a random GM fourfold $X\subset Y$ containing 
the image $S=\overline{\psi(V)}\subset Y$ of the surface $V$ 
(the lack of error messages means that $X$ is smooth).
{\footnotesize
\begin{Verbatim}[commandchars=&!$]
&colore!darkorange$!i3 :$ psi = &colore!airforceblue$!rationalMap$(&colore!airforceblue$!ideal$ B,&colore!airforceblue$!Dominant$=>2);
&colore!circOut$!o3 =$ &colore!output$!RationalMap (cubic rational map from PP^5 to 5-dimensional subvariety of PP^8)$
&colore!darkorange$!i4 :$ S = psi(&colore!airforceblue$!ideal$ V);
&colore!darkorange$!i5 :$ X = &colore!bleudefrance$!specialGushelMukaiFourfold$ S;
&colore!circOut$!o5 :$ &colore!output$!ProjectiveVariety, GM fourfold containing a surface of degree 7 and sectional genus 0$
\end{Verbatim}
} \noindent 
The discriminant of $X$ can be obtained with one of the two commands:
{\footnotesize
\begin{Verbatim}[commandchars=&!$]
&colore!darkorange$!i6 :$ &colore!bleudefrance$!discriminant$ X
&colore!circOut$!o6 =$ &colore!output$!20$
&colore!darkorange$!i7 :$ &colore!bleudefrance$!describe$ X
&colore!circOut$!o7 =$ &colore!output$!Special Gushel-Mukai fourfold of discriminant 20$
&colore!output$!     containing a surface in PP^8 of degree 7 and sectional genus 0$
&colore!output$!     cut out by 21 hypersurfaces of degree 2$
&colore!output$!     and with class in G(1,4) given by 4*s_(3,1)+3*s_(2,2)$
\end{Verbatim}
} \noindent 
The following code tells us that 
$S$ corresponds 
to a smooth point of an irreducible family $\mathcal S$ of dimension $21$ 
of surfaces in $Y$, and that the family of GM fourfolds containing some surface of $\mathcal S$
has codimension $3$ in the space of all GM fourfolds (we omit a few lines of output).
{\footnotesize
\begin{Verbatim}[commandchars=&!$]
&colore!darkorange$!i8 :$ &colore!bleudefrance$!parameterCount$ X
&colore!output$!-- h^1(N_{S,Y}) = 0$
&colore!output$!-- h^0(N_{S,Y}) = 21$
&colore!output$!-- h^1(O_S(2)) = 0 and h^0(I_{S,Y}(2)) = 16 = h^0(O_Y(2)) - \chi(O_S(2))$
&colore!output$!-- h^0(N_{S,X}) = 0$
&colore!circOut$!o8 =$ &colore!output$!(3,(16,21,0))$
\end{Verbatim}
} \noindent 
Finally, the following code reveals a congruence  
of $3$-secant conics of the surface $S$ inside $Y$. 
Indeed, the rational map $\phi:Y\dashrightarrow\PP^{15}$ 
defined by the quadrics through $S$ is birational onto a fivefold 
$Z\subset\PP^{15}$, and through the general point $\phi(p)$ of $Z$ 
there pass $5$ lines contained in $Z$. Of these $5$ lines, only one comes from a $3$-secant conic to 
$S$ passing through $p$.
 Note that this implies the rationality of $X$,
 a fact already known for the GM fourfolds of discriminant $20$ (see \cite{HoffSta}).
{\footnotesize
\begin{Verbatim}[commandchars=&!$]
&colore!darkorange$!i9 :$ &colore!bleudefrance$!detectCongruence$ X
&colore!output$!-- phi: quadratic rational map from Y to PP^15, Z = phi(Y)$
&colore!output$!-- number lines contained in Z and passing through the point phi(p): 5$
&colore!output$!-- number 1-secant lines to S passing through p: 4$
&colore!output$!-- number 3-secant conics to S passing through p: 1$
\end{Verbatim}


\begin{table}[htbp]
\renewcommand{\arraystretch}{1.300}
\centering
\tabcolsep=0.3pt 
\begin{adjustbox}{width=\textwidth}
\begin{tabular}{|c|c|c|c|c|c|c|c|c|c|}
\hline
\rowcolor{gray!5.0}
 & {\begin{tabular}{c} Surface \\ $S\subset Y=\mathbb{G}(1,4)\cap\PP^8$ \end{tabular}} & $K_S^2$ & {\begin{tabular}{c} Class in \\ $\mathbb{G}(1,4)$\end{tabular}} & {\begin{tabular}{c} Codim \\ in $\mathcal M_{4}^{GM}$ \end{tabular}} & \begin{tabular}{c} Image \\ in $\mathcal D$\end{tabular} & $h^0(\mathcal I_{S/Y}(2))$ & $h^0(N_{S/Y})$ & $h^0(N_{S/X})$ &
  \scriptsize{\begin{tabular}{c} Curves of degree $e$ in $Y$\\ passing though a general \\ point of $Y$ and that  are \\ $(2e-1)$-secant to $S$ \\ for $e=1,2,3$  \end{tabular}} \\
\hline \hline 
1 & {\begin{tabular}{c} $\tau$-quadric surface  \cite{DIM} \end{tabular}} & $8$  & $\sigma_{3,1}+\sigma_{2,2}$ & $1$ & $\mathcal D_{10}'$ & $31$ & $8$ & $0$ & $1$, $0$, $0$  \\
\hline
2 & {\begin{tabular}{c} Quartic scroll \end{tabular}}  &  $8$ & $3 \sigma_{3,1} + \sigma_{2,2}$ & $1$ & $\mathcal D_{10}'$ & $25$ & $14$ & $0$ & $3$, $0$, $0$   \\
\hline
3 &  {\begin{tabular}{c} K3 surface of degree \\ $14$ and genus $8$  \cite{ProcRS} \end{tabular}}  & $0$  & $9 \sigma_{3,1}+ 5 \sigma_{2,2}$ & $1$ & $\mathcal D_{10}'$ & $10$ & $39$ & $10$ & $9$, $8$, $1$   \\
\hline
4 & {\begin{tabular}{c} Quintic del Pezzo \\ surface  \cite{Roth1949} \end{tabular}}  & $5$  & $3 \sigma_{3,1}+2 \sigma_{2,2}$ & $1$ & $\mathcal D_{10}''$ & $24$ & $18$ & $3$ & $3$, $0$, $0$  \\
\hline 
5 & {\begin{tabular}{c} Rational surface of \\ degree $9$ and genus $3$ \end{tabular}}  & $1$  & $5 \sigma_{3,1}+4 \sigma_{2,2}$ & $1$ & $\mathcal D_{10}''$ & $16$ & $28$ & $5$ & $5$, $0$, $0$  \\
\hline
6 & {\begin{tabular}{c} $\sigma$-plane  \cite{Roth1949} \end{tabular}}  & $9$  & $\sigma_{3,1}$ & $2$ & $\mathcal D_{10}''$ & $34$ & $4$ & $0$ & $1$, $0$, $0$ \\
\hline
7 & {\begin{tabular}{c} Cubic scroll  \cite{DIM} \end{tabular}}  & $8$  & $2 \sigma_{3,1} + \sigma_{2,2}$ & $1$ & $\mathcal D_{12}$ & $28$ & $11$ & $0$ & $2$, $0$, $0$  \\
\hline 
8 & {\begin{tabular}{c} Rational surface of \\ degree $7$ and genus $2$ \end{tabular}}  &  $3$ & $4 \sigma_{3,1}+3 \sigma_{2,2}$ & $1$ & $\mathcal D_{12}$ & $20$ & $23$ & $4$ & $4$, $0$, $0$ \\
\hline
9 & {\begin{tabular}{c} $\rho$-plane  \cite{Roth1949} \end{tabular}}  & $9$  & $\sigma_{2,2}$ & $3$ & $\mathcal D_{12}$ & $34$ & $3$ & $0$ & $0$, $0$, $0$ \\
\hline 
10 & {\begin{tabular}{c} Rational surface of \\ degree $10$ and genus $4$ \end{tabular}}  & $0$ & $6 \sigma_{3,1}+4 \sigma_{2,2}$ & $1$ & $\mathcal D_{16}$ & $15$ & $29$ & $5$ & $6$, $0$, $0$  \\
\hline 
11 & {\begin{tabular}{c} Rational surface of \\ degree $10$ and genus $3$ \end{tabular}}  & $2$ & $6 \sigma_{3,1}+4 \sigma_{2,2}$ & $1$ & $\mathcal D_{16}$ & $13$ & $29$ & $3$ & $6$, $2$, $0$ \\
\hline 
12 & {\begin{tabular}{c} K3 surface of degree \\ $14$ and genus $8$ \end{tabular}}  & $0$ & $8 \sigma_{3,1}+ 6 \sigma_{2,2}$ & $1$ & $\mathcal D_{16}$ & $10$ & $38$ & $9$ & $8$, $8$, $0$  \\
\hline
13 & {\begin{tabular}{c} Rational surface of \\ degree $12$ and genus $5$ \end{tabular}}  & $-1$ & $7 \sigma_{3,1}+ 5 \sigma_{2,2}$ & $1$ & $\mathcal D_{18}'$ & $11$ & $32$ & $4$ & $7$, $5$, $0$ \\
\hline 
14 & {\begin{tabular}{c} Rational surface of \\ degree $8$ and genus $2$ \end{tabular}}  & $4$ & $5 \sigma_{3,1}+3 \sigma_{2,2}$ & $\leq 2$ & $\mathcal D_{18}'$ & $17$ & $24$ & $3$ & $5$, $0$, $0$ \\
\hline
15 & {\begin{tabular}{c} Rational surface of \\ degree $9$ and genus $3$ \end{tabular}}  & $2$ & $5 \sigma_{3,1} + 4 \sigma_{2,2}$ & $1$ & $\mathcal D_{18}''$ & $16$ & $26$ & $3$ & $5$, $0$, $0$ \\
\hline 
16 & {\begin{tabular}{c} Rational surface of \\ degree $11$ and genus $5$ \end{tabular}}  & $-1$ & $7 \sigma_{3,1}+4 \sigma_{2,2}$ & $1$ & $\mathcal D_{18}''$ & $14$ & $30$ & $5$ & $7$, $0$, $0$ \\
\hline 
17 & {\begin{tabular}{c} Rational surface of \\ degree $9$ and genus $2$ \\ \cite{HoffSta} \end{tabular}}  & $5$ & $6 \sigma_{3,1}+3 \sigma_{2,2}$ & \begin{tabular}{c} $\leq 2$ \\\tiny{(really $1$)} \end{tabular} & $\mathcal D_{20}$ & $14$ & $25$ & \begin{tabular}{c} $1$ \\ \tiny{(but $0$ when } \\ \tiny{$S$ is general)} \end{tabular} & $6$, $1$, $0$ \\
\hline 
18 & {\begin{tabular}{c} Septic scroll \end{tabular}}  & $8$ & $4 \sigma_{3,1} + 3 \sigma_{2,2}$ & $3$ & $\mathcal D_{20}$ & $16$ & $21$ & $0$ & $4$, $1$, $0$ \\
\hline 
19 & {\begin{tabular}{c} Rational surface of \\ degree $10$ and genus $3$ \end{tabular}}  & $3$ & $6 \sigma_{3,1}+4 \sigma_{2,2}$ & $1$ & $\mathcal D_{24}$ & $13$ & $27$ & $1$ & $6$, $2$, $0$\\
\hline
20 & {\begin{tabular}{c} Veronese surface \end{tabular}}  & $9$ & $2 \sigma_{3,1}+2 \sigma_{2,2}$ & $4$ & $\mathcal D_{24}$ & $25$ & $11$ & $0$ & $2$, $0$, $0$ \\
\hline 
21 & {\begin{tabular}{c} Rational surface of \\ degree $12$ and genus $5$ \end{tabular}}  & $0$ & $7 \sigma_{3,1}+ 5 \sigma_{2,2}$ & $1$ & $\mathcal D_{26}'$ & $11$ & $30$ & $2$ & $7$, $4$, $0$ \\
\hline 
\end{tabular}
\end{adjustbox}
 \caption{Families of  GM fourfolds described as the closure
 of the  locus of smooth quadric hypersurfaces in a del Pezzo fivefold $Y$ 
 containing some smooth surface $S\subset Y$ 
 varying in an irreducible component of 
 $\mathrm{Hilb}_Y^{\chi(\mathcal{O}_S(t))}$.}
 \label{tabella1} 
\end{table}

\begin{table}[htbp]
\renewcommand{\arraystretch}{1.019}
\centering
\tabcolsep=0.3pt 
\begin{adjustbox}{width=\textwidth}
\begin{tabular}{|c|c|c|c|c|c|c|c|c|}
\hline
\rowcolor{gray!5.0}
 & $B\subset\PP^5$ & $V\subset\PP^5$ & $C\subseteq B\cap V$  & $B\cap V\setminus C$ & $V\simeq \psi_B(V)$ & $\dim\{B:B\supset C\}$ & $\dim\{V:V\supset C\}$ & \begin{tabular}{c}\begin{tabular}{l} $\dim\{(C,B,V):$ \end{tabular}\\ \begin{tabular}{r} $C \subseteq B\cap V\}$\end{tabular} \end{tabular}\\
\hline \hline 
1 & \begin{tabular}{c} cubic scroll \\ surface \end{tabular} & \begin{tabular}{c} quadric surface \end{tabular} & \begin{tabular}{c} irreducible conic \\ curve \end{tabular} & $\emptyset$ & yes & $11$ & $6$ & $31$  \\
\hline 
2 & \begin{tabular}{c} cubic scroll \\ surface \end{tabular} & \begin{tabular}{c} quadric surface \end{tabular} & \begin{tabular}{c} line of the \\ ruling \end{tabular} & $\emptyset$  & yes & $16$ & $10$ & $34$ \\
\hline 
3 & \begin{tabular}{c} cubic scroll \\ surface \end{tabular} & \begin{tabular}{c} K3 surface of degree $8$ \\ and genus $5$ \end{tabular} & \begin{tabular}{c} union of $3$ lines \\ of the ruling with \\ the directrix line \end{tabular} & $\emptyset$  & yes & $3$ & $28$ & $57$ \\
\hline 
4 & \begin{tabular}{c} cubic scroll \\ surface \end{tabular} & \begin{tabular}{c} del Pezzo surface \\ of degree $4$  \end{tabular} & \begin{tabular}{c} twisted cubic \\ curve  \end{tabular} & $\emptyset$  & no & $7$ & $13$ & $40$ \\
\hline 
5 & \begin{tabular}{c} quartic scroll \\ threefold \end{tabular} & \begin{tabular}{c} cubic scroll \\ surface  \end{tabular} & \begin{tabular}{c} irreducible conic \\ curve  \end{tabular} & $7$ points  & no & $22$ & $11$ & $47$ \\
\hline 
6 & \begin{tabular}{c} cubic scroll \\ surface \end{tabular} & \begin{tabular}{c} plane  \end{tabular} & \begin{tabular}{c} line of the \\ruling \end{tabular} & $\emptyset$  & yes & $16$ & $3$ & $27$ \\
\hline
7 & \begin{tabular}{c} cubic scroll \\ surface \end{tabular} & \begin{tabular}{c} cubic scroll \\ surface  \end{tabular} & \begin{tabular}{c} twisted cubic \\curve  \end{tabular} & $\emptyset$  & yes & $7$ & $7$ & $34$ \\
\hline 
8 & \begin{tabular}{c} quartic scroll \\ threefold \end{tabular} & \begin{tabular}{c} quadric surface  \end{tabular} & \begin{tabular}{c} line \end{tabular} & $5$ points  & no & $25$ & $10$ & $43$ \\
\hline 
9 & \begin{tabular}{c} cubic scroll \\ surface \end{tabular} & \begin{tabular}{c} plane  \end{tabular} & \begin{tabular}{c} directrix line \end{tabular} & $\emptyset$  & yes & $15$ & $3$ & $26$ \\
\hline 
10 & \begin{tabular}{c} quartic scroll \\ threefold \end{tabular} & \begin{tabular}{c} del Pezzo surface \\ of degree $4$ \end{tabular} & \begin{tabular}{c} rational normal \\ quartic curve \end{tabular} & $4$ points  & no & $14$ & $8$ & $48$ \\
\hline
11 & \begin{tabular}{c} quartic scroll \\ threefold \end{tabular} & \begin{tabular}{c} quartic scroll \\ surface \end{tabular} & \begin{tabular}{c} rational normal \\ quartic curve \end{tabular} & $6$ points & no & $14$ & $8$ & $48$ \\
\hline 
12 & \begin{tabular}{c} cubic scroll \\ surface \end{tabular} & \begin{tabular}{c} K3 surface of degree $8$ \\ and genus $5$ \end{tabular} & \begin{tabular}{c} rational normal \\ quartic curve \end{tabular} & $\emptyset$  & yes & $2$ & $28$ & $56$  \\
\hline 
13 & \begin{tabular}{c} cubic scroll \\ surface \end{tabular} & \begin{tabular}{c} rational surface of \\ degree $7$ and genus $3$, \\ the image of the plane \\via the linear system \\of quintic curves with \\ nine simple base points \\and one triple point \end{tabular} & \begin{tabular}{c} rational normal \\ quartic curve, the \\ image of a general \\ conic passing through \\ three of the nine base \\ points and through \\the triple point  \end{tabular} & $\emptyset$  & yes & $2$ & $22$ & $50$  \\
\hline 
14 & \begin{tabular}{c} cubic scroll \\ surface \end{tabular} & \begin{tabular}{c} del Pezzo surface \\ of degree $4$ \end{tabular} & \begin{tabular}{c} irreducible conic \\ curve \end{tabular} & $\emptyset$  & yes & $11$ & $18$ & $43$ \\
\hline 
15 & \begin{tabular}{c} cubic scroll \\ surface \end{tabular} & \begin{tabular}{c} rational surface of \\ degree $6$ and genus $2$, \\ the image of the plane \\ via the linear system \\ of quartic curves with \\ six simple base points \\and one double point \end{tabular} & \begin{tabular}{c} rational normal \\ quartic curve, the \\image of a general \\ conic passing through \\ two of the six base \\ points and through \\the double point \end{tabular} & $\emptyset$  & yes & $2$ & $17$ & $45$ \\
\hline 
16 & \begin{tabular}{c} cubic scroll \\ surface \end{tabular} & \begin{tabular}{c} rational surface of \\ degree $6$ and genus $3$,\\ the image of the plane \\ via the linear system\\of quartic curves with \\ten simple base points \end{tabular} & \begin{tabular}{c} twisted cubic \\ curve, the image \\ of a general conic \\passing through five  \\ of the ten base points \end{tabular} & $\emptyset$ & yes & $7$ & $21$ & $48$ \\
\hline 
17 & \begin{tabular}{c} cubic scroll \\ surface \end{tabular} & \begin{tabular}{c} del Pezzo surface \\ of degree $5$ \end{tabular} & \begin{tabular}{c} twisted cubic \\ curve \end{tabular} & $\emptyset$  & yes & $7$ & $17$ & $44$ \\
\hline 
18 & \begin{tabular}{c} quartic scroll \\ threefold \end{tabular} & \begin{tabular}{c} quartic scroll \\ surface \end{tabular} & \begin{tabular}{c} rational normal \\ quintic curve \end{tabular} & $3$ points & no & $11$ & $2$ & $45$ \\
\hline 
19 & \begin{tabular}{c} quartic scroll \\ threefold \end{tabular} & \begin{tabular}{c} Veronese surface \end{tabular} & \begin{tabular}{c} rational normal \\ quartic curve \end{tabular} & $6$ points  & no & $14$ & $6$ & $46$ \\
\hline 
20 & \begin{tabular}{c} cubic scroll \\ surface \end{tabular} & \begin{tabular}{c} Veronese surface \end{tabular} & \begin{tabular}{c} rational normal \\ quartic curve \end{tabular} & $\emptyset$  & yes & $2$ & $6$ & $34$ \\
\hline 
21 & \begin{tabular}{c} cubic scroll \\ surface \end{tabular} & \begin{tabular}{c} rational surface of \\ degree $7$ and genus $3$, \\ the image of the plane \\via the linear system \\of quartic curves with \\ nine simple base points \end{tabular} & \begin{tabular}{c} rational normal \\ quartic curve, the \\ image of a general \\ conic passing through \\ four of the nine base \\ points  \end{tabular} & $\emptyset$  & yes & $2$ & $20$ & $48$ \\
\hline 
\end{tabular}
\end{adjustbox}
 \caption{Construction of surfaces  $S\subset Y$ as in Table~\ref{tabella1}, obtained 
 as images 
 of 
 surfaces $V\subset\PP^5$ via 
  the 
  map $\psi_B:\PP^5\dashrightarrow Y$ associated to $B$.}
 \label{tabella2} 
\end{table}

\clearpage 

\begin{table}[htbp]
\renewcommand{\arraystretch}{1.300}
\centering
\tabcolsep=1.1pt 
\begin{adjustbox}{width=\textwidth}
\begin{tabular}{|c|c|c|c|c|c|c|c|c|c|c|c|c|c|c|c|c|c|c|c|c|}
\hline
\rowcolor{gray!5.0}
$1$ & $2$ & $3$ & $4$ & $5$ & $6$ & $7$ & $8$ & $9$ & $10$ & $11$ & $12$ & $13$ & $14$ & $15$ & $16$ & $17$ & $18$ & $19$ & $20$ & $21$ \\
\hline 
$14s$ & $7m$ & $18m$ & $12s$ & $9.4h$ & $25s$ & $16s$ & $5h$ & $2s$ & $2h$ & $2.5h$ & $11.2h$ & $5.5h$ & $2h$ & $20m$ & $37m$ & $10m$ & $7m$ & $17m$ & $8m$ & $30m$ \\
\hline
\end{tabular}
\end{adjustbox}
\caption{Run-times to execute the function \texttt{parameterCount}
in a code like the one in Section~\ref{sec computations},
for all the examples in Tables~\ref{tabella1} and \ref{tabella2},
on a laptop with an Intel Core i7-6500U processor, 16 GB of RAM.}
\end{table}


\providecommand{\bysame}{\leavevmode\hbox to3em{\hrulefill}\thinspace}
\providecommand{\MR}{\relax\ifhmode\unskip\space\fi MR }
\providecommand{\MRhref}[2]{%
  \href{http://www.ams.org/mathscinet-getitem?mr=#1}{#2}
}
\providecommand{\href}[2]{#2}

\end{document}